\documentclass[10pt,twoside,reqno]{amsart}
\usepackage{amssymb}
\usepackage{multirow}
\calclayout

\numberwithin{equation}{section}
\makeatletter

\renewcommand{\@secnumfont}{\bfseries}

\renewcommand{\section}{\@startsection{section}{1}%
  {0mm}{.7\linespacing\@plus\linespacing}{.5\linespacing}
  {\normalfont\bfseries\centering}}

\newcommand{\bibsection}{\@startsection{section}{1}%
  {0mm}{.7\linespacing\@plus\linespacing}{.5\linespacing}
  {\normalfont\scshape\centering}}

\renewcommand{\@biblabel}[1]{#1.}

\newtheorem{thm}{\bf Theorem}[section]

\begin{document}

\vspace{1.3cm}

\title{Degenerate Daehee polynomials of the second kind}

\author{Taekyun Kim}
\address{Department of Mathematics, Kwangwoon University, Seoul 139-701, Republic
    of Korea}
\email{tkkim@kw.ac.kr}

\author{Dae San Kim}
\address{Department of Mathematics, Sogang University, Seoul 121-742, Republic
    of Korea}
\email{dskim@sogang.ac.kr}

\begin{abstract}
In this paper, we consider the degenerate Daehee numbers and
polynomials of the second kind which are different from the
previously introduced degenerate Daehee numbers and polynomials.
We investigate some properties of these numbers and polynomials.
In addition, we give some new identities and relations between the
Daehee polynomials of the second kind and Carlitz's degenerate
Bernoulli polynomials.
\end{abstract}

\maketitle

\markboth{\centerline{\scriptsize Degenerate Daehee polynomials of the second kind}}
{\centerline{\scriptsize T. Kim, D. S. Kim}}

\bigskip
\medskip

\section{\bf Introduction}

As is well known, the Bernoulli polynomials are defined by the generating function
\begin{equation} \begin{split} \label{01}
\frac{t}{e^{t}-1}e^{xt}=\sum_{n=0}^\infty B_{n}(x)\frac{t^n}{n!},\quad (\textnormal{see} \,\, [1,20]).
\end{split} \end{equation}
When $x=0$, $B_{n}=B_{n}(0)$ are called the Bernoulli numbers.

In [3], L.Carlitz considered the degenerate Bernoulli polynomials given by
\begin{equation} \begin{split} \label{02}
\frac{t}{(1+\lambda t)^{\frac{1}{\lambda}}-1}(1+\lambda t)^{\frac{x}{\lambda}}=\sum_{n=0}^\infty \beta_{n,\lambda}(x)\frac{t^n}{n!},\,\,\,\,(\lambda \in \mathbb{R}).
\end{split} \end{equation}
When $x=0$, $\beta_{n,\lambda}=\beta_{n,\lambda}(0)$ are called the degenerate Bernoulli numbers.

The falling factorial sequence is given by
\begin{equation} \begin{split} \label{03}
(x)_{0}=1,\,\,\,\,(x)_{n}=x(x-1)\cdots(x-n+1),\,\,\,\,(n \geq 1).
\end{split} \end{equation}

The stirling numbers of the first kind are defined by
\begin{equation} \begin{split} \label{04}
(x)_{n}=\sum_{l=0}^{n}S_{1}(n,l)x^{l},\,\,\,\,(n\geq0),\quad (\textnormal{see} \,\, [6,7,8,13,20]).
\end{split} \end{equation}

It is well known that the stirling numbers of the second kind are defined as
\begin{equation*} \begin{split}
x^{n}=\sum_{l=0}^{n}S_{2}(n,l)(x)_{l},\,\,\,\,(n \geq 0),\quad (\textnormal{see} \,\, [6,12,20]).
\end{split} \end{equation*}

From \eqref{01} and \eqref{02}, we note that
\begin{equation} \begin{split} \label{05}
\sum_{n=0}^\infty \lim_{\lambda \rightarrow 0}\beta_{n,\lambda}(x)\frac{t^n}{n!}&=\lim_{\lambda \rightarrow 0}\frac{t}{(1+\lambda t)^{\frac{1}{\lambda}}-1}(1+\lambda t)^{\frac{x}{\lambda}}\\
&=\frac{t}{e^t-1}e^{xt}=\sum_{n=0}^\infty B_{n}(x)\frac{t^n}{n!}.
\end{split} \end{equation}

Thus, by \eqref{05}, we get
\begin{equation} \begin{split} \label{06}
\lim_{\lambda \rightarrow 0}\beta_{n,\lambda}(x)=B_{n}(x),\,\,\,\,(n\geq 0).
\end{split} \end{equation}

The $\lambda$-analogue of falling factorial sequence is defined by L.Carlitz as follows:
\begin{equation} \begin{split} \label{07}
&(x)_{0,\lambda}=1,\,\,\,\, (x)_{n,\lambda}=x(x-\lambda)\cdots(x-(n-1)\lambda),\,\,(n\geq 1),\quad (\textnormal{see} \,\, [3]).
\end{split} \end{equation}

By\eqref{02}, we get
\begin{equation} \begin{split} \label{08}
\sum_{n=0}^\infty \beta_{n,\lambda}(x)\frac{t^n}{n!}&=\bigg(\sum_{m=0}^\infty \lambda^{-m}B_{m}(x)\frac{\Big(\log(1+\lambda t)\Big)^m}{m!}\bigg)\bigg(\frac{\lambda t}{\log(1+\lambda t)}\bigg)\\
&=\bigg(\sum_{m=0}^\infty \lambda^{-m}B_{m}(x)\sum_{k=m}^\infty S_{1}(k,m)\frac{\lambda^{k}t^{k}}{k!}\bigg)\bigg(\sum_{l=0}^\infty b_{l}\frac{\lambda^{l}t^{l}}{l!}\bigg)\\
&=\bigg(\sum_{k=0}^\infty\Big(\sum_{m=0}^{k}\lambda^{k-m}B_{m}(x)S_{1}(k,m)\Big)\frac{t^k}{k!}\bigg)\bigg(\sum_{l=0}^\infty b_{l}\frac{\lambda^{l}t^{l}}{l!}\bigg)\\
&=\sum_{n=0}^\infty\bigg(\sum_{k=0}^{n}\sum_{m=0}^{k}{n \choose k}\lambda^{n-m}B_{m}(x)S_{1}(k,m)b_{n-k}\bigg)\frac{t^n}{n!},
\end{split} \end{equation}

Here $b_{n},(n\geq 0)$, are the Bernoulli numbers of the second kind given by the generating function
\begin{equation} \begin{split} \label{09}
\frac{t}{\log(1+t)}=\sum_{n=0}^\infty b_{n}\frac{t^n}{n!},\quad (\textnormal{see} \,\, [1,20]).
\end{split} \end{equation}

From \eqref{08}, we note that
\begin{equation} \begin{split} \label{10}
\beta_{n,\lambda}(x)=\sum_{k=0}^{n}\sum_{m=0}^{k}{n \choose k}\lambda^{n-m}B_{m}(x)S_{1}(k,m)b_{n-k}.
\end{split} \end{equation}

The Daehee polynomials are defined by
\begin{equation} \begin{split} \label{11}
\frac{\log(1+t)}{t}(1+t)^{x}=\sum_{n=0}^\infty D_{n}(x)\frac{t^n}{n!},\quad (\textnormal{see} \,\, [10,21-23]).
\end{split} \end{equation}
When $x=0$, $D_{n}=D_{n}(0)$ are called the Daehee numbers.

From \eqref{01} and \eqref{11}, we have
\begin{equation} \begin{split} \label{12}
\sum_{n=0}^\infty B_{n}(x)\frac{t^n}{n!}&=\sum_{m=0}^\infty D_{m}(x)\frac{1}{m!}(e^t-1)^m\\
&=\sum_{m=0}^\infty D_{m}(x)\sum_{n=m}^\infty S_2(n,m)\frac{t^n}{n!}\\
&=\sum_{n=0}^\infty\bigg(\sum_{m=0}^{n}S_{2}(n,m)D_{m}(x)\bigg)\frac{t^n}{n!}.
\end{split} \end{equation}

Comparing the coefficients on both sides of \eqref{12}, we have
\begin{equation} \begin{split} \label{13}
B_{n}(x)=\sum_{m=0}^{n}S_{2}(n,m)D_{m}(x),\,\,(n\geq 0),\quad (\textnormal{see} \,\, [4,6,12]).
\end{split} \end{equation}

The degenerate Daehee numbers are given by the generating function
\begin{equation} \begin{split} \label{14}
\frac{\lambda \log\big(1+\frac{1}{\lambda}\log(1+\lambda t)\big)}{\log(1+\lambda t)}=\sum_{n=0}^\infty \tilde{D}_{n,\lambda}\frac{t^n}{n!},\quad (\textnormal{see} \,\, [14]).
\end{split} \end{equation}

Recently, many authors have studied the Daehee numbers and polynomials and the degenerate Daehee numbers and polynomials (see [2-19,21-23]).\\
In this paper, we consider the degenerate Daehee numbers and polynomials of the second kind which are different from the previous introduced degenerate Daehee numbers and polynomials. We investigate some properties of these numbers and polynomials. In addition, we give some new identities and relations between the Daehee polynomials of the second kind and Carlitz's degenerate Bernoulli polynomials.

\section{Degenerate Daehee polynomials}

For $\lambda \in \mathbb{R}$, we consider the degenerate Daehee polynomials of the second kind given by the generating function
\begin{equation} \begin{split} \label{15}
\frac{\log(1+t)}{\big(1+\lambda \log(1+t)\big)^{\frac{1}{\lambda}}-1}\big(1+\lambda \log(1+t)\big)^{\frac{x}{\lambda}}=\sum_{n=0}^\infty D_{n,\lambda}(x)\frac{t^n}{n!}.
\end{split} \end{equation}
When $x=0$, $D_{n,\lambda}=D_{n,\lambda}(0)$ are called the degenerate Daehee numbers of the second kind.

Note that
\begin{equation} \begin{split} \label{16}
\sum_{n=0}^\infty \lim_{\lambda \rightarrow 0} D_{n,\lambda}(x)\frac{t^n}{n!}&=\lim_{\lambda \rightarrow 0}\frac{\log(1+t)}{\big(1+\lambda \log(1+t)\big)^{\frac{1}{\lambda}}-1}\big(1+\lambda \log(1+t)\big)^{\frac{x}{\lambda}}\\
&=\frac{\log(1+t)}{t}(1+t)^{x}=\sum_{n=0}^\infty D_{n}(x)\frac{t^n}{n!}.
\end{split} \end{equation}

Thus, by \eqref{16}, we get
\begin{equation*} \begin{split}
\lim_{\lambda \rightarrow 0}D_{n,\lambda}(x)=D_{n}(x),\,\,\,\,(n \geq 0).
\end{split} \end{equation*}

From \eqref{02}, we note that
\begin{equation} \begin{split} \label{17}
&\frac{\log(1+t)}{\big(1+\lambda \log(1+t)\big)^{\frac{1}{\lambda}}-1}\big(1+\lambda \log(1+t \big)^{\frac{x}{\lambda}}\\
&=\sum_{m=0}^\infty \beta_{m,\lambda}(x)\frac{1}{m!}\big(log(1+t)\big)^{m}\\
&=\sum_{m=0}^\infty \beta_{m,\lambda}(x)\sum_{n=m}^\infty S_{1}(n,m)\frac{t^n}{n!}\\
&=\sum_{n=0}^\infty\big(\sum_{m=0}^{n}\beta_{m,\lambda}(x)S_{1}(n,m)\big)\frac{t^n}{n!}.
\end{split} \end{equation}

Therefore, by \eqref{15} and \eqref{17}, we obtain the following theorem.
\begin{thm} For $n \geq 0$, we have
\begin{equation*} \begin{split}
D_{n,\lambda}(x)=\sum_{m=0}^{n}\beta_{m,\lambda}(x)S_{1}(n,m).
\end{split} \end{equation*}
\end{thm}

We observe that
\begin{equation}\begin{split}\label{18}
&\frac{\log(1+t)}{\big(1+\lambda \log(1+t)\big)^{\frac{1}{\lambda}}-1}\big(1+\lambda \log(1+t)\big)^{\frac{x}{\lambda}}\\
&=\frac{\log(1+t)}{\big(1+\lambda \log(1+t)\big)^{\frac{1}{\lambda}}-1} \sum_{l=0}^\infty \left( \frac{x}{\lambda} \right)_l \lambda^l \frac{1}{l!} \big( \log(1+t)\big)^\lambda\\
&=\left( \sum_{m=0}^\infty D_{m,\lambda} \frac{t^m}{m!} \right) \left( \sum_{l=0}^\infty (x)_{l,\lambda} \sum_{k=l}^\infty S_1(k,l) \frac{t^k}{k!} \right)\\
&=\left( \sum_{m=0}^\infty D_{m,\lambda} \frac{t^m}{m!} \right) \left( \sum_{k=0}^\infty \left( \sum_{l=0}^k (x)_{l,k}  S_1(k,l) \right) \frac{t^k}{k!}\right)\\
&=\sum_{n=0}^\infty \left( \sum_{k=0}^n \sum_{l=0}^k {n \choose k} (x)_{l,\lambda} S_1(k,l) D_{n-k,\lambda} \right) \frac{t^n}{n!}.
\end{split}\end{equation}

Therefore, by \eqref{15} and \eqref{18}, we obtain the following theorem.

\begin{thm} For $n \geq 0$, we have
\begin{equation*}\begin{split}
D_{n,\lambda}(x) = \sum_{k=0}^n \sum_{l=0}^k {n \choose k} (x)_{l,\lambda} S_1(k,l) D_{n-k,\lambda}.
\end{split}\end{equation*}
\end{thm}

By replacing $t$ by $e^t-1$ in \eqref{15}, we get
\begin{equation}\begin{split}\label{19}
\sum_{m=0}^\infty D_{m,\lambda}(x) \frac{1}{m!} (e^t-1)^m &=
\frac{t}{(1+\lambda t)^{\frac{1}{\lambda}}-1}(1+\lambda t)^{\frac{x}{\lambda}}\\
&= \sum_{n=0}^\infty \beta_{n,\lambda}(x) \frac{t^n}{n!}.
\end{split}\end{equation}

On the other hand,
\begin{equation}\begin{split}\label{20}
\sum_{m=0}^\infty D_{m,\lambda}(x) \frac{1}{m!} (e^t-1)^m &= \sum_{m=0}^\infty D_{m,\lambda }(x) \sum_{n=m}^\infty   S_2(n,m) \frac{t^n}{n!}\\
&= \sum_{n=0}^\infty \left( \sum_{m=0}^n D_{m,\lambda }(x) S_2(n,m) \right) \frac{t^n}{n!}.
\end{split}\end{equation}
Therefore, by \eqref{19} and \eqref{20}, we obtain the following theorem.

\begin{thm} For $n \geq 0$, we have
\begin{equation*}\begin{split}
 \beta_{n,\lambda}(x) = \sum_{m=0}^n D_{m,\lambda }(x) S_2(n,m).
\end{split}\end{equation*}
\end{thm}

By \eqref{15}, we get

\begin{equation}\begin{split}\label{21}
\log(1+t) &= \sum_{n=0}^\infty D_{n,\lambda }\frac{t^n}{n!}\left( \big( 1+\lambda \log(1+t) \big)^{\frac{1}{\lambda }}-1 \right)\\
&= \frac{\log(1+t)}{\big(1+\lambda \log(1+t)\big)^{\frac{1}{\lambda}}-1}\big( 1+\lambda \log(1+t) \big)^{\frac{1}{\lambda }}-\sum_{n=0}^\infty D_{n,\lambda }\frac{t^n}{n!}\\
&=\sum_{n=0}^{\infty}(D_{n,\lambda }(1) - D_{n,\lambda })\frac{t^n}{n!}.
\end{split}\end{equation}

On the other hand
\begin{equation}\begin{split}\label{22}
\log(1+t) = \sum_{n=1}^\infty \frac{(-1)^{n-1}}{n}t^n.
\end{split}\end{equation}
Therefore, by \eqref{21} and \eqref{22}, we obtain the following theorem.

\begin{thm}
For $n \geq 0$, we have
\begin{equation*}\begin{split}
D_{n,\lambda }(1) - D_{n,\lambda } = \begin{cases}
0,&\text{if}\,\,n=0,\\
(-1)^{n-1} (n-1)!,& \text{if}\,\, n \geq 1.
\end{cases}
\end{split}\end{equation*}
\end{thm}

From \eqref{15}, we have

\begin{equation}\begin{split}\label{23}
&\frac{\log(1+t)}{\big(1+\lambda \log(1+t)\big)^{\frac{1}{\lambda }}-1}
\big(1+\lambda \log(1+t)\big)^{\frac{x}{\lambda }}\\
&=
\frac{\log(1+t)}{\big(1+\lambda \log(1+t)\big)^{\frac{d}{\lambda }}-1} \sum_{a=0}^{d-1} \big(1+\lambda \log(1+t)\big)^{\frac{a+x}{\lambda }}\\
&= \frac{1}{d} \left( \frac{d \log(1+t)}{\big(1+ \frac{\lambda }{d}(d \log(1+t) \big)^{\frac{d}{\lambda }}-1 } \right)\sum_{a=0}^{d-1} \big( 1+ \frac{\lambda }{d} (d \log(1+t) \big)^{\frac{d}{\lambda } \left( \frac{a+x}{d} \right)}\\
&= \frac{1}{d} \sum_{a=0}^{d-1} \sum_{m=0}^\infty \beta_{m,\frac{\lambda }{d}}\left( \frac{a+x}{d} \right) \frac{1}{m!} \big(d\log(1+t) \big)^m\\
&= \frac{1}{d} \sum_{a=0}^{d-1} \sum_{m=0}^\infty \beta_{m,\frac{\lambda }{d}}\left( \frac{a+x}{d} \right)  \sum_{n=m}^\infty d^m S_1(n,m) \frac{t^n}{n!}\\
&= \sum_{n=0}^\infty \left\{ \sum_{m=0}^n \sum_{a=0}^{d-1} d^{m-1}\beta_{m,\frac{\lambda }{d}}\left( \frac{a+x}{d} \right) S_1(n,m) \right\} \frac{t^n}{n!},
\end{split}\end{equation}
where $d \in \mathbb{N}$ and $n \geq 0$.

By \eqref{15} and \eqref{23}, we get the following theorem.

\begin{thm} For $d \in \mathbb{N}$ and $n \geq 0$, we have
\begin{equation*}\begin{split}
D_{n,\lambda }(x) = \sum_{m=0}^n d^{m-1} S_1(n,m)\sum_{a=0}^{d-1} \beta_{m,\frac{\lambda }{d}}\left( \frac{a+x}{d} \right).
\end{split}\end{equation*}
\end{thm}

On the one hand, we have

\begin{equation}\begin{split}\label{24}
&\log(1+t) \sum_{l=0}^{n-1} \big(1+\lambda \log(1+t)\big)^{\frac{l}{\lambda }}\\
&= \frac{\log(1+t)}{\big(1+\lambda \log(1+t)\big)^{\frac{1}{\lambda }}-1}
\big(1+\lambda \log(1+t)\big)^{\frac{n}{\lambda }}- \frac{\log(1+t)}{\big(1+\lambda \log(1+t)\big)^{\frac{1}{\lambda }}-1}\\
&= \sum_{m=0}^\infty \left\{ D_{m,\lambda }(n) - D_{m,\lambda } \right\} \frac{t^m}{m!} = t \sum_{m=0}^\infty \left\{ \frac{D_{m+1,\lambda}(n) - D_{m+1,\lambda } }{m+1} \right\} \frac{t^m}{m!}.
\end{split}\end{equation}

On the other hand, we have

\begin{equation}\begin{split}\label{25}
&\log(1+t) \sum_{l=0}^{n-1} \big(1+\lambda \log(1+t)\big)^{\frac{l}{\lambda }} = t \left( \frac{\log(1+t)}{t} \right) \sum_{l=0}^{n-1}  \big( 1+\lambda \log(1+t) \big)^{\frac{l}{\lambda }}\\
&= t \left( \sum_{p=0}^\infty D_p \frac{t^p}{p!} \right) \left( \sum_{j=0}^\infty \left( \sum_{k=0}^j (l)_{k,\lambda } S_1(j,k) \right) \frac{t^j}{j!} \right) \\
&= t \sum_{m=0}^\infty \left( \sum_{j=0}^m \sum_{k=0}^j {m \choose j} (l)_{k,\lambda } S_1(j,k) D_{m-j} \right) \frac{t^m}{m!}.
\end{split}\end{equation}

Therefore, by \eqref{24} and \eqref{25}, we obtain the following theorem.

\begin{thm} For $n \in \mathbb{N}$, $m \geq 0$, we have
\begin{equation*}\begin{split}
\frac{1}{m+1} \left( D_{m+1,\lambda}(n) - D_{m+1,\lambda } \right)  = \sum_{j=0}^m \sum_{k=0}^j {m \choose j} (l)_{k,\lambda } S_1(j,k) D_{m-j} .
\end{split}\end{equation*}
\end{thm}

For $r \in \mathbb{N}$, we define the higher-order degenerate Daehee polynomials of the second kind given by the generating function

\begin{equation}\begin{split}\label{26}
&\left(\frac{\log(1+t)}{\big(1+\lambda \log(1+t)\big)^{\frac{1}{\lambda }}-1}\right)^r \big(1+\lambda \log(1+t)\big)^{\frac{x}{\lambda }}= \sum_{n=0}^\infty D_{n,\lambda }^{(r)}(x) \frac{t^n}{n!}.
\end{split}\end{equation}

When $x=0$, $D_{n,\lambda }^{(r)} = D_{n,\lambda }^{(r)}(0)$ are called the higher-order degenerate Daehee numbers of the second kind.

From \eqref{26}, we note that

\begin{equation*}\begin{split}
\sum_{n=0}^\infty \lim_{\lambda \rightarrow 0} D_{n,\lambda }^{(r)}(x) \frac{t^n}{n!} &= \lim_{\lambda \rightarrow 0} \left(\frac{\log(1+t)}{\big(1+\lambda \log(1+t)\big)^{\frac{1}{\lambda }}-1}\right)^r \big(1+\lambda \log(1+t)\big)^{\frac{x}{\lambda }}\\
&= \left( \frac{\log(1+t)}{t} \right)^r (1+t)^x = \sum_{n=0}^{\infty}D_n^{(r)}(x) \frac{t^n}{n!},
\end{split}\end{equation*}

\noindent where $D_n^{(r)}(x)$ are called the higher-order Daehee polynomials.

As is well known, the higher-order degenerate Bernoulli polynomials are considered by L. Carlitz as follows:

\begin{equation}\begin{split}\label{27}
\left( \frac{t}{(1+\lambda t)^{\frac{1}{\lambda }}-1}\right)^r (1+\lambda t)^{\frac{x}{\lambda }} = \sum_{n=0}^\infty \beta_{n,\lambda }^{(r)}(x) \frac{t^n}{n!}.
\end{split}\end{equation}
Note that $\lim_{\lambda  \rightarrow 0} \beta_{n,\lambda }^{(r)}(x) = B_n^{(r)}(x)$, $(n \geq 0)$, where $B_n^{(r)}(x)$ are the higher-order Bernoulli polynomials.

From \eqref{26}, we note that

\begin{equation}\begin{split}\label{28}
&\left(\frac{\log(1+t)}{\big(1+\lambda \log(1+t)\big)^{\frac{1}{\lambda }}-1}\right)^r \big(1+\lambda \log(1+t)\big)^{\frac{x}{\lambda }}\\
&= \sum_{m=0}^\infty \beta_{m,\lambda }^{(r)}(x) \frac{1}{m!} \big( \log(1+t)\big)^m = \sum_{m=0}^\infty \beta_{m,\lambda }^{(r)}(x) \sum_{n=m}^\infty S_1(n,m) \frac{t^n}{n!}\\
&= \sum_{n=0}^\infty \left( \sum_{m=0}^n \beta_{m,\lambda }^{(r)} S_1(n,m) \right) \frac{t^n}{n!}.
\end{split}\end{equation}

Thus, by \eqref{27} and \eqref{28}, we obtain the following theorem.

\begin{thm} For $n \geq 0$, we have
\begin{equation*}\begin{split}
D_{n,\lambda }^{(r)}(x) = \sum_{m=0}^n \beta_{m,\lambda }^{(r)} S_1(n,m).
\end{split}\end{equation*}
\end{thm}

By replacing $t$ by $e^t-1$ in \eqref{26}, we get

\begin{equation}\begin{split}\label{29}
\sum_{m=0}^\infty D_{m,\lambda }^{(r)}(x) \frac{1}{m!} (e^t-1)^m &= \left( \frac{t}{(1+\lambda t)^{\frac{1}{\lambda }}-1} \right)(1+\lambda t)^{\frac{x}{\lambda }}\\
&=\sum_{n=0}^\infty \beta_{n,\lambda }^{(r)}(x) \frac{t^n}{n!}.
\end{split}\end{equation}

On the other hand,

\begin{equation}\begin{split}\label{30}
\sum_{m=0}^\infty D_{m,\lambda }^{(r)}(x) \frac{1}{m!} (e^t-1)^m &= \sum_{m=0}^\infty D_{m,\lambda }^{(r)}(x) \sum_{n=m}^\infty S_2(n,m) \frac{t^n}{n!}\\
&= \sum_{n=0}^\infty \left( \sum_{m=0}^n D_{m,\lambda }^{(r)}(x) S_2(n,m) \right) \frac{t^n}{n!}.
\end{split}\end{equation}

Therefore, by \eqref{29} and \eqref{30}, we obtain the following theorem.

\begin{thm}
For $n \geq 0$, we have
\begin{equation*}\begin{split}
\beta_{m,\lambda }^{(r)}(x) =  \sum_{m=0}^n D_{m,\lambda }^{(r)}(x) S_2(n,m).
\end{split}\end{equation*}
\end{thm}

For $r,k \in \mathbb{N}$, with $r>k$, by \eqref{26}, we get

\begin{equation}\begin{split}\label{31}
&\left(\frac{\log(1+t)}{\big(1+\lambda \log(1+t)\big)^{\frac{1}{\lambda }}-1}\right)^r \big(1+\lambda \log(1+t)\big)^{\frac{x}{\lambda }}\\
&=\left(\frac{\log(1+t)}{\big(1+\lambda \log(1+t)\big)^{\frac{1}{\lambda }}-1}\right)^{r-k} \left(\frac{\log(1+t)}{\big(1+\lambda \log(1+t)\big)^{\frac{1}{\lambda }}-1}\right)^{k}\big(1+\lambda \log(1+t)\big)^{\frac{x}{\lambda }}\\
&=\left( \sum_{l=0}^\infty D_{l,\lambda }^{(r-k)} \frac{t^l}{l!}  \right)
\left( \sum_{m=0}^\infty D_{m,\lambda }^{(k)}(x) \frac{t^m}{m!}  \right)\\
&= \sum_{n=0}^\infty \left( \sum_{l=0}^n {n \choose l} D_{l,\lambda }^{(r-k)} D_{n-l,\lambda }^{(k)}(x) \right) \frac{t^n}{n!}.
\end{split}\end{equation}

Therefore, by \eqref{26} and \eqref{31}, we obtain the following theorem.

\begin{thm}
For $r,k \in \mathbb{N}$, with $r>k$, we have
\begin{equation*}\begin{split}
D_{n,\lambda }^{(r)}(x) = \sum_{l=0}^n {n \choose l} D_{l,\lambda }^{(r-k)} D_{n-l,\lambda }^{(k)}(x),\,\,(n \geq 0).
\end{split}\end{equation*}
\end{thm}

It is well known that

\begin{equation}\begin{split}\label{32}
\left( \frac{t}{\log(1+t)} \right)^k (1+t)^{x-1} = \sum_{n=0}^\infty B_n^{(n-k+1)}(x) \frac{t^n}{n!},\,\,(k \in \mathbb{Z}),
\end{split}\end{equation}
where $B_n^{(\alpha)}(x)$ are called the higher-order Bernoulli polynomials which are given by the generating function

\begin{equation*}\begin{split}
\left( \frac{t}{e^t-1} \right)^{\alpha} e^{xt} = \sum_{n=0}^\infty B_n^{(\alpha)} (x) \frac{t^n}{n!}.
\end{split}\end{equation*}
Thus, by \eqref{32}, we get

\begin{equation}\begin{split}\label{33}
&\left(\frac{\log(1+t)}{\big(1+\lambda \log(1+t)\big)^{\frac{1}{\lambda }}-1}\right)^r \big(1+\lambda \log(1+t)\big)^{\frac{x}{\lambda }}\\
&=\left( \frac{\big(1+\lambda \log(1+t)\big)^{\frac{1}{\lambda }}-1}{\frac{1}{\lambda}\log \big(1+\lambda \log(1+t)\big)} \right)^{-r} \big(1+\lambda \log(1+t)\big)^{\frac{x}{\lambda }} \left( \frac{\lambda \log(1+t)}{\log(1+\lambda \log(1+t))} \right)^r \\
&= \left(\sum_{m=0}^\infty B_m^{(m+r+1)}(x+1) \frac{1}{m!} \left( \big(1+\lambda \log(1+t)\big)^{\frac{1}{\lambda }}-1 \right)^m \right)\\&\quad \times \left( \frac{\lambda \log(1+t)}{\log(1+\lambda \log(1+t))} \right)^r.
\end{split}\end{equation}

From \eqref{32}, we note that

\begin{equation}\begin{split}\label{34}
\left( \frac{\lambda \log(1+t)}{\log(1+\lambda \log(1+t))} \right)^r &=
\sum_{j=0}^\infty B_j^{(j-r+1)}(1) \frac{1}{j!}\lambda ^j \big(\log(1+t)\big)^j\\&= \sum_{l=0}^\infty \left( \sum_{j=0}^l B_j^{(j-r+1)}(1) \lambda ^j S_1(l,j) \right) \frac{t^l}{l!}.
\end{split}\end{equation}

As is known, the degenerate Stirling numbers of the second kind are defined by the generating function

\begin{equation}\begin{split}\label{35}
\frac{1}{m!} \big( (1+\lambda t)^{\frac{1}{\lambda }}-1 \big)^m = \sum_{n=m}^\infty S_{2,\lambda }(n,m) \frac{t^n}{n!},
\end{split}\end{equation}
where $m \in \mathbb{N}$ with $m \geq 0$, (see [7]).

Note that $\lim_{\lambda  \rightarrow 0} S_{2,\lambda }(n,m) = S_2(n,m)$, $(n,m \geq 0)$. Also, we note that

\begin{equation}\begin{split}\label{36}
&\sum_{m=0}^\infty B_m^{(m+r+1)}(x+1) \frac{1}{m!} \left( \big(1+\lambda \log(1+t)\big)^{\frac{1}{\lambda }}-1\right)^m \\
&= \sum_{m=0}^\infty B_m^{(m+r+1)}(x+1) \sum_{k=m}^\infty S_{2,\lambda }(k,m) \frac{1}{k!} \big( \log(1+t)\big)^k \\
&=\sum_{k=0}^\infty \left( \sum_{m=0}^k B_m^{(m+r+1)}(x+1)S_{2,\lambda }(k,m) \right) \sum_{p=k}^\infty S_1(p,k) \frac{t^p}{p!}\\
&= \sum_{p=0}^\infty \left\{ \sum_{k=0}^p \sum_{m=0}^k B_m^{(m+r+1)}(x+1) S_{2,\lambda }(k,m) S_1(p,k) \right\} \frac{t^p}{p!}.
\end{split}\end{equation}

From \eqref{33}, \eqref{34}, and \eqref{36}, we have

\begin{equation}\begin{split}\label{37}
&\left(\frac{\log(1+t)}{\big(1+\lambda \log(1+t)\big)^{\frac{1}{\lambda }}-1}\right)^r \big(1+\lambda \log(1+t)\big)^{\frac{x}{\lambda }}\\
&= \left(\sum_{p=0}^\infty \left\{ \sum_{k=0}^p \sum_{m=0}^k B_m^{(m+r+1)}(x+1) S_{2,\lambda }(k,m) S_1(p,k) \right\} \frac{t^p}{p!} \right)\\&\quad \times \left( \sum_{l=0}^\infty \left( \sum_{j=0}^l B_j^{(j-r+1)}(1) \lambda ^j S_1(l,j) \right) \frac{t^l}{l!} \right)\\
&= \sum_{n=0}^\infty \Bigg( \sum_{p=0}^n \sum_{k=0}^p \sum_{m=0}^k \sum_{j=0}^{n-p} {n \choose p} B_m^{(m+r+1)}(x+1) B_j^{(j-r+1)}(1) \lambda ^{j} S_{2,\lambda }(k,m) S_1(p,k)\\
&\quad \times  S_1(n-p,j) \Bigg) \frac{t^n}{n!}
\end{split}\end{equation}

Therefore, by \eqref{26} and \eqref{37}, we get the following result.

\begin{thm} For $n \geq 0$, we have
\begin{equation*}\begin{split}
D_{n,\lambda }^{(r)}(x) &= \sum_{p=0}^n \sum_{k=0}^p \sum_{m=0}^k \sum_{j=0}^{n-p} {n \choose p} B_m^{(m+r+1)}(x+1) B_j^{(j-r+1)}(1) \lambda ^{j} S_{2,\lambda }(k,m) \\
&\quad \times S_1(p,k) S_1(n-p,j).
\end{split}\end{equation*}
\end{thm}

From \eqref{26}, we note that

\begin{equation}\begin{split}\label{38}
&\left(\frac{\log(1+t)}{\big(1+\lambda \log(1+t)\big)^{\frac{1}{\lambda }}-1}\right)^r \big(1+\lambda \log(1+t)\big)^{\frac{x}{\lambda }}\\
&= \left( \sum_{l=0}^\infty D_{l,\lambda }^{(r)}\frac{t^l}{l!} \right) \left(
\sum_{k=0}^\infty \left( \sum_{m=0}^k (x)_{m,\lambda }S_1(k,m) \right) \frac{t^k}{k!} \right)\\
&= \sum_{n=0}^\infty  \left( \sum_{k=0}^n \sum_{m=0}^k {n \choose k} (x)_{m,\lambda } S_1(k,m) D_{n-k,\lambda}^{(r)} \right) \frac{t^n}{n!}.
\end{split}\end{equation}

Thus, by \eqref{26} and \eqref{38}, we get the next theorem.

\begin{thm} For $n \geq 0$, we have
\begin{equation*}\begin{split}
D_{n,\lambda }^{(r)}(x) = \sum_{k=0}^n \sum_{m=0}^k {n \choose k} (x)_{m,\lambda } S_1(k,m) D_{n-k,\lambda}^{(r)}.
\end{split}\end{equation*}
\end{thm}

Now, we observe that

\begin{equation}\begin{split}\label{39}
\sum_{n=0}^\infty D_{n,\lambda }^{(r)}(x+y) \frac{t^n}{n!} &= \left(\frac{\log(1+t)}{\big(1+\lambda \log(1+t)\big)^{\frac{1}{\lambda }}-1}\right)^r \big(1+\lambda \log(1+t)\big)^{\frac{x+y}{\lambda }}\\
&= \left( \sum_{l=0}^\infty D_{l,\lambda }^{(r)}(x) \frac{t^l}{l!} \right) \left( \sum_{k=0}^\infty \left( \sum_{m=0}^k (y)_{m,\lambda } S_1(k,m) \right) \frac{t^k}{k!} \right)\\
&= \sum_{n=0}^\infty \left( \sum_{k=0}^n \sum_{m=0}^k {n \choose k} D_{n-k}^{(r)}(x) (y)_{m,\lambda } S_1(k,m) \right) \frac{t^n}{n!}.
\end{split}\end{equation}

Thus, by \eqref{39}, we get

\begin{equation}\begin{split}\label{40}
D_{n,\lambda }^{(r)}(x+y) = \sum_{k=0}^n \sum_{m=0}^k {n \choose k} D_{n-k}^{(r)}(x) (y)_{m,\lambda } S_1(k,m).
\end{split}\end{equation}

From \eqref{32}, we note that

\begin{equation}\begin{split}\label{41}
&\left( \frac{\big(1+\lambda \log(1+t)\big)^{\frac{1}{\lambda }}-1}{\log(1+t)}\right)^r = \left( \frac{t}{\log(1+t)} \right)^r \frac{r!}{t^r}\frac{1}{r!} \left( \big(1+\lambda \log(1+t)\big)^{\frac{1}{\lambda }}-1\right)^r\\
&=\left( \sum_{m=0}^\infty B_m^{(m-r+1)}(1) \frac{t^m}{m!} \right) \frac{r!}{t^r} \left( \sum_{l=r}^\infty S_{2,\lambda }(l,r) \frac{1}{l!} \big( \log(1+t)\big)^l \right)\\
&=\left( \sum_{m=0}^\infty B_m^{(m-r+1)}(1) \frac{t^m}{m!} \right) \frac{r!}{t^r} \left( \sum_{l=0}^\infty S_{2,\lambda }(l+r,r) \frac{1}{(l+r)!} \big( \log(1+t)\big)^{l+r} \right) \\
&=\left( \sum_{m=0}^\infty B_m^{(m-r+1)}(1) \frac{t^m}{m!} \right) \frac{r!}{t^r} \left( \sum_{k=r}^\infty \left( \sum_{l=0}^{k-r} S_{2,\lambda }(l+r,r) S_1(k,l+r) \right) \frac{t^k}{k!} \right)\\
&=\left( \sum_{m=0}^\infty B_m^{(m-r+1)}(1) \frac{t^m}{m!} \right) \left( \sum_{k=0}^\infty \left( \sum_{l=0}^k \frac{S_{2,\lambda }(l+r,r)S_1(k+r,l+r)}{{k+r \choose k}}\right) \frac{t^k}{k!} \right)\\
&= \sum_{n=0}^\infty \left\{ \sum_{k=0}^n \sum_{l=0}^k \frac{{n \choose k}}{{k+r \choose k}} S_{2,\lambda }(l+r,r) S_1(k+r,l+r) B_{n-k}^{(n-k-r+1)}(1) \right\} \frac{t^n}{n!}.
\end{split}\end{equation}

On the other hand,

\begin{equation}\begin{split}\label{42}
\left( \frac{\big(1+\lambda \log(1+t)\big)^{\frac{1}{\lambda }}-1}{\log(1+t)}\right)^r = \sum_{n=0}^\infty D_{n,\lambda }^{(-r)} \frac{t^n}{n!}.
\end{split}\end{equation}

Therefore, by \eqref{41} and \eqref{42}, we obtain the following theorem.

\begin{thm}
For $n \geq 0$, $r \in \mathbb{N}$, we have
\begin{equation*}\begin{split}
 D_{n,\lambda }^{(-r)} = \sum_{k=0}^n \sum_{l=0}^k \frac{{n \choose k}}{{k+r \choose k}} S_{2,\lambda }(l+r,r) S_1(k+r,l+r) B_{n-k}^{(n-k-r+1)}(1).
\end{split}\end{equation*}
\end{thm}

\end{document}